\documentstyle[twoside]{article}

\newtheorem{thm}{\indent{\sc Theorem}}[section]
\newtheorem{defn}[thm]{\indent{\sc Definition}} 
\newtheorem{lem}[thm]{\indent{\sc Lemma}} 
\newtheorem{por}[thm]{\indent{\sc Porism}}
\newtheorem{prop}[thm]{\indent{\sc Proposition}} 
\newtheorem{cor}[thm]{\indent{\sc Corollary}}
\newtheorem{ex}{\indent {\sc Example}}[section]
\newtheorem{prob}[thm]{\indent{\sc Open Problem}\index{Problem{,} Open}}

\newcommand{\ssect}{\subsection}
\newcommand{\sssect}{\subsubsection}

\newcommand{\eqref}[1]{equation~(\ref{#1})}
\newcommand{\eref}[1]{\eqref{#1}}

\newcommand{~}{\nolinebreak[3] }

\newcommand{\th}{\raisebox{0.6ex}{th}}

\newcommand{\seper}{\nopagebreak \begin{center}
	\underline{\hspace{2in}}
	\end{center}}

\newcommand{\emp}[1]{{\em #1}\index{#1}}
\newcommand{\bold}[1]{{\em (#1)}\index{#1}}
\newcommand{\bld}[1]{{\em (#1)}}

 \newcommand{\proof}[1]{\proo{#1}$\Box $}
\newcommand{\proo}[1]{{\em Proof: }#1}
\renewcommand{\Box}{ \hspace{1cm}{\rule{1.2ex}{2ex}}}

\newcommand{\D}{{{\rm D}}}

\newcommand{\SB}[1]{\mbox{{\bf \scriptsize #1}}}

\newcommand{\digress}[1]{ \begin{quotation}
	 #1 {\em End of Digression.} \end{quotation}}

\renewcommand{\hat}[1]{\widehat{#1}}

\begin{document}
\include{title}

\begin{abstract}
We define the Artinian and Noetherian algebra which consist of formal series
involving exponents which are not necessarily integers. All of the usual
operations are defined here and characterized. As an application, we compute
the algebra of symmetric functions with nonnegative real exponents. The
applications to logarithmic series and the Umbral calculus are deferred to
another paper.\cite{ch4}
\seper
{\bf Les Series \`{a} Exposants Quelconques}

On d\'{e}finit ici les alg\`{e}bres Artinienne et Noetherienne comme
\'{e}tant des alg\`{e}bres constitu\'{e}es des s\'{e}ries formelles
\`{a} exposants pas n\'{e}cessairement entiers. On definit sur ces
alg\`{e}bres toutes les op\'{e}rations classiques et on les
caracterise. Comme exemple d'exploitation de cette th\'{e}orie, on
s'interesse \`{a} alg\`{e}bre de fonctions sym\'{e}triques \`{a}
exponsants r\`{e}els en nonn\'{e}gatifs. Une autre publication
\cite{ch4} est 
consacr\'{e}e aux applications aux series logarithmiques et au calcul
ombral. 
\end{abstract}

\begin{center}
{\it Dedicated to\\
H\'el\`ene}
\end{center}

\tableofcontents

\section{Introduction}

The simplest type of series is the polynomial. However, it is common to
occasionally consider other more general series formal power series, Laurent
series, inverse formal power series and so on. Nevertheless, little
consideration has been given to series whose exponents are not necessarily
integers. 
We would like to derive a ``continuous'' analog of Laurent or inverse Laurent
series in which exponents are chosen from the real numbers (or any other poset)
instead of merely from the integers.
Here for any choice of coefficients and exponents we define two sets
of formal series: {\em Artinian} and {\em Noetherian}. 

These series are not merely of academic interest, since the Noetherian series
in infinitely many variables 
represent asymptotic expansions of real functions with respect to the ladder of
comparison 
$$x^{a}(\log x)^{b}(\log \log x)^{c}\cdots $$
where $a,b,c,\ldots $
are real numbers, and the Artinian series in the derivative represent the
logarithmic analog of shift-invariant operators. We defer this application 
to another article \cite{ch4}. 

However, as typical applications of this theory, we
will compute the Artinian and Noetherian algebras of 
symmetric functions with nonnegative real exponents, and derive the Lagrange
inversion formula for Artinian and Noetherian series.

For these applications, it must be proven
that one can manipulate formal power series with real exponents as easily as
one does polynomials. To that end, we define {\em Artinian series} in which for
all real numbers 
$a$ there are finitely many terms of degree $b$ with $b<a$, and {\em Noetherian
series} which have the dual condition.

They are equipped with a topology and operations which make them a topological
algebra over the complex numbers. Moreover, both types of series actually form
a field. 

Even further operations can be defined. We define $f(x)$ raised to any real
exponent $a$. This definition involves the choice of any arbitrary integer $n$, so we
write $f(x)^{a;n}$. The operation is then characterized in terms of its
algebraic properties. We further define the composition of one series with
another $f(g;n)$ and characterize this. This composition is not always
associative; however it is associative in the cases relevant to \cite{ch4}.

\ssect{The Artinian Algebra and Noetherian Algebra} 

\begin{defn}\label{artseq}
\bold{Artinian Sequence}
Let $R$ be a poset (usually the real
numbers). Define an  \emp{Artinian poset}\index{Artinian poset} to be any
subposet $S\subseteq R$ such that 
for all $a\in S$ there are only finitely many $b\leq a$ such
that $b\in S$. That is, all principal ideals of the
subposet $S$ are finite.

Let $K$ be an additive group (usually the complex numbers).
We say a sequence $(c_{a})_{a\in R}$ of group elements indexed by
the poset is {\em Artinian} if its
support is Artinian. That is, if  for all $a\in 
R$ there are only finitely many $b\leq a$ such that $c_{b}\neq 0$.

Dually, we say $S$ is a {\em Noetherian poset}\index{Noetherian Poset} if for
all $a\in S$   
there are only finitely many $b\geq a$ such that $b\in S$.
That is, all principal filters of the subposet $S$ are finite.
We  say a sequence is {\em
Noetherian}\index{Noetherian Sequence} if its support is 
Noetherian. That is, if for all $a\in R$
there are only finitely many $b\geq a$ such that $c_{b}\neq 0$.
\end{defn}

We note a few trivial observations. Nonempty Artinian
(resp.\ Noetherian) posets have  minimal (resp.\ maximal)
elements, so  nonzero Artinian (resp.\ Noetherian) 
sequences have  lowest nonzero terms. 

A poset is both Artinian and Noetherian if and only if it is finite, so
a sequence is both Artinian 
and Noetherian if and only if it has finite support.

\begin{defn} \bold{Artinian Algebra}\label{ArtAlg}
Define the {\em Artinian algebra} 
$K(x)^{R}$ to be the set of all {\em formal} sums $f(x)=\sum_{a\in
\SB{R}}c_{a}x^{a}$ such that the $c_{a}$ form an Artinian
sequence. $f(x)$ is called a {\em Artinian series.}\index{Artinian Series}

Define the {\em degree}\index{Degree}
of $f(x)$ to be 
$$ \deg(f(x))=\left\{ \begin{array}{ll}
\max \{a:c_{a}\neq 0 \}&\mbox{if $f(x)\neq 0$, and}\\*[0.1in]
-\infty &\mbox{if $f(x)=0$.}
\end{array} \right. $$

Similarly, the {\em Noetherian algebra}\index{Noetherian Algebra}  $K(x)_{R}$
is the collection of all  
formal sums $f(x)=\sum _{a\in R}c_{a}x^{a}$ where the $c_{a}$ form a Noetherian
sequence. 
$f(x)$ is called a {\em Noetherian series}\index{Noetherian Series.} Define its
{\em  degree}\index{Degree}  to be
$$ \deg(f(x))=\left\{ \begin{array}{ll}
\min \{a:c_{a}\neq 0 \}&\mbox{if $f(x)\neq 0$, and}\\*[0.1in]
+\infty &\mbox{if $f(x)=0$.}
\end{array} \right. $$

In either case, we  denote the {\em coefficient}\index{Coefficient}
of $x^{a}$ in $f(x)$ by $[x^{a}]f(x)$. For $f(x)\neq 0$, the coefficient
$\left[x^{\deg(f(x))}\right]f(x)$ 
is called the {\em leading coefficient}\index{Leading Coefficient} of
$f(x)$, and
$\left(\left[x^{\deg(f(x))}\right]f(x)\right)x^{\deg(f(x))}$
is called the  
{\em leading term}\index{Leading Term} of $f(x)$.  A series of degree one is
called a {\em delta series.}\index{Delta Series}
\end{defn}

\begin{prop} \label{ArtAlgP} Suppose $K$ is a ring, and $R$ is an
ordered monoid. Then the Artinian Algebra and the Noetherian
Algebra are $K$-algebras when they are equipped with the
operations: 
\begin{eqnarray}
\label{op1}
 \left( \sum _{a}c_{a}x^{a}\right)+
\left(\sum_{a}d_{a}x^{a}\right) & =& \sum
_{a}(c_{a}+d_{a})x^{a}\\*
\label{op2}
\left( \sum _{a}c_{a}x^{a}\right)\left(\sum_{a}d_{a}x^{a}\right)
&=& \sum _{c}\left(\sum _{a+b=c}c_{a}d_{b} \right)x^{c}\\*
\label{op3}
z\left(\sum _{a}c_{a}x^{a} \right) &=& \sum _{a}(zc_{a})x^{a}.
\end{eqnarray}
\end{prop}

\proof{It must be shown that all these operations are well
defined, and that they obey the axioms of $K$-algebras.  For
brevity, we only show that multiplication is well defined in
the Noetherian Algebra.
To do this, we must first show that the coefficients of a
product are well defined, and then show that the
coefficients form a Noetherian sequence. 

Let $e_{c}=\sum _{a+b=c}c_{a}d_{b}$. To show that $e_{c}$ is well defined, let
$a_{0}$ be the maximum $a$ such that $c_{a}\neq 0$. Now, 
there are only finitely many $b\geq c-a_{0}$ such that $d_{b}\neq 0$.
Thus, $e_{c}$ is well defined.

Let $b_{0}$ be the maximum $b$ such that $d_{b}\neq 0$. Chose some
$c'\geq c$ such that $e_{c'}\neq 0$. Then there are some $a$ and
$b$ which sum to $c'$ and such that $c_{a}\neq 0$ and $d_{b}\neq
0$. Moreover, $a>c-b_{0}$ and $b>c-a_{0}$, so there is only
finitely many such $a$ and $b$. Thus, there are only finitely many
such $c'$. Hence, the sequence $e_{c}$ is Noetherian.}

\begin{cor} Suppose $K$ is an integral domain, and $R$ is an
ordered monoid. Then the Artinian algebra and Noetherian
algebra are integral domains.$\Box $ \end{cor}

\begin{prop} \label{Inverses} Suppose $K$ is a field, and $R$ is
an Archimedean group.  Then the Artinian Algebra and the
Noetherian Algebra are fields.  \end{prop}

\proof{Given Proposition~\ref{ArtAlgP}, it suffices only to
calculate inverses. We only treat the Noetherian case, since
the Artinian case is similar.

Let $f(x)$ be a nonzero Noetherian series of degree $a$ with coefficients given
by  $c_{b}=[x^{b}]f(x)$. We calculate the coefficients of its inverse. Let
$d_{-a}=c_{a}^{-1},$ and for $b<-a$, 
let 
$$ d_{b} =-c_{a}^{-1}\sum _{b< e\leq -a}d_{b}c_{a+b-e}. $$
 $d_{b}$ is well defined,
since there are finitely many nonzero terms being summed over in
each step of the recursion. We claim that $d_{b}$ is Noetherian.
If so, it is easy to check that $\sum _{b} d_{b}x^{b}$ is
the multiplicative inverse of $f(x)$. 

Suppose $d_{b}$ is not Noetherian. Then there is a sequence
$b_{1}<b_{2}<\cdots <-a$ such that $d_{b_{i}}\neq 0$ for
all $i$. Since the reals are Archimedean, $b_{i}-b_{i+1}$ must
tend to zero (in the order topology).  However $d_{b}\neq 0$ for
$b<-a$ only if $d_{a+b-e}\neq 0$ where $e$ is
one of the finitely many $e\geq a-b$ such that $c_{e}\neq
0$. The set of differences among the various $e$ is finite so it
has a lower bound. Contradiction.}

Note that Theorem~\ref{Charact2} provides an easier proof of
a more sweeping result under slightly stronger conditions.
Conversely, the following Porism strips the requirements for
the existence of inverses to the bone.

\begin{por} \label{poor}
Suppose $K$ is a ring, and $R$ is an Archimedean monoid.
Then a nonzero Artinian or Noetherian series has a multiplicative
inverse if and only if its degree has an additive inverse and its
leading coefficient has a multiplicative inverse.$\Box $ 
\end{por}

Some of the most common types are series are Noetherian or Artinian series:
\begin{ex}
 $K(x)^{{\bf N}}$ is merely
the ring of polynomials in the variable $x$. The only invertible
series are invertible constants.  
\end{ex}

\begin{ex}
$K(x)^{{\bf Z}}$ (where
$K$ is a field) is the field of Laurent\index{Laurent
series} series in the variable $x$.
\end{ex}

\begin{ex}
$K(x)_{{\bf N}}$ is the ring of formal power series in
the variable $x$. The only invertible series are ones with an
invertible constant term.  
\end{ex}

We wish to define series in several variables. However, we must be careful,
since the Artinian algebra $\left(K(x)^{R} \right)(y)^{R}$ is not necessarily
equal to the Artinian algebra $\left(K(y)^{R} \right)(x)^{R}$. For example,
\begin{eqnarray*} 
\sum _{n\geq 0} x^{n}y^{-n} &\in
&\left(K(x)^{{\bf Z}} \right)(y)^{{\bf Z}}\\ 
&\not\in &
\left(K(y)^{{\bf Z}} \right)(x)^{{\bf Z}}.  
\end{eqnarray*}

\begin{defn}
\bold{Multivariate Noetherian Algebra} \label{MNA}
Let ${\cal R}=\{ R_{1},\ldots ,R_{n} \}$ be a collection of
posets and let 
${\cal X}=\{x_{1},\ldots ,x_{n}\}$ be a collection of variables. Given a
group $K$ we recursively define the {\em multivariate
Noetherian algebra}\index{Noetherian Algebra,Multivariate}
in the variables ${\cal X}$ indexed by ${\cal R}$ by the recursion
$$ K(x_{1;R_{1}},\ldots ,x_{n;R_{n}})= \left\{ \begin{array}{ll}
\left(K(x_{n})_{R_{n}} \right)(x_{1;R_{1}},\ldots
,x_{n-1;R_{n-1}}) &\mbox{for $n>0$, and}\\
\\
K&\mbox{for $n=0$.}
\end{array} \right. $$

 Let
${\cal R}=\{R_{1},R_{2},\ldots  \}$ be an infinite
collection of posets and let 
${\cal X}=\{x_{1},x_{2},\ldots \}$ be an infinite collection
of variables. The {\em infinite multivariate Noetherian
algebra}\index{Noetherian Algebra}
in the variables ${\cal X}$ indexed by ${\cal R}$ is defined
to be the direct limit ({\em ie:} union) of the finite
multivariate Noetherian algebras
$ K(x_{1;R_{1}},x_{2;R_{2}},\ldots ) = \bigcup_{n\geq 0}
K(x_{1;R_{1}},\ldots ,x_{n;R_{n}})$ 

When $R_{n}=R$ for all $n$, we write $K(x_{1},\ldots
,x_{n})_{R}$ for the finite multivariate Noetherian algebra,
and $K(x_{1},x_{2},\ldots )_{R}$ for the infinite
multivariate Noetherian algebra.
\end{defn}

We chose this particular definition for the multivariate Noetherian algebra
because in \cite{ch4} the continuous Logarithmic Algebra ${\cal I}= {\bf
C}(x,\log x,\log \log x ,\ldots )_{\bf R}$ defined in terms of it contains only
asymptotic expansions as $x$ tends towards infinity whereas not all members
of ${\bf C}(\log \log x,\log x,x)_{\bf R}$ represent any sort of asymptotic
expansion. 

An Artinian algebra of several variables $K(x_{1},\ldots
,x_{n})^{R}$ or $K(x_{1},x_{2},\ldots )^{R}$ could be defined
similarly,\index{Artinian Algebra, Multivariate} for there is a strong duality
between the Noetherian and Artinian algebra: 

\begin{prop} \label{Duality} When $K$ is a ring, and $R$ is an ordered
group, there is a canonical isomorphism $\iota
:K(x)_{R}\rightarrow K(x)^{R}$ given by $\iota :x^{a}\mapsto
x^{-a}\Box $ \end{prop}

Under the topology to be defined in the next section, this isomorphism is
actually a topological algebra isomorphism.

\ssect{Topology} 
\index{Topology}

We now give an alternate definition of the
Artinian and Noetherian algebras as the topological completion of
a simpler algebra. This algebra called the {\em finite algebra} is
actually the intersection of the Artinian and Noetherian
algebras.

\begin{defn} \bold{Finite Algebra}\label{ArtTop}
Suppose $K$ is a topological ring, and $R$ is a
ordered monoid.  Define the {\em  finite algebra} $K[x;R]$ to be the set
of all {\em finite} sums $\sum _{a\in R}c_{a}x^{a}$. That is
$c_{a}\in K$ is nonzero only for finitely many values of $a$.
Members of the finite algebra are called {\em finite series}\index{Finite
Series}.

The Noetherian (resp.\ Artinian) topology for $K[x;R]$ is the
finest topology such that 
\begin{enumerate} 
\item $c_{n}x^{a_{n}}$
is a Cauchy sequence whenever 
\begin{enumerate}
\item $a_{n}$ is a Cauchy sequence in the
order topology of $R$, and 
\item  $c_{n}$ is a Cauchy sequence in the
topology of $K$.  
\end{enumerate}
\item $c_{n}x^{a_{n}}$ converges to zero
whenever $a_{n}$ decreases (resp.\ increases) without bound.  \item
Finite sums of the above two items are also Cauchy sequences.
\end{enumerate} \end{defn}

\begin{prop} 
The finite algebra is a topological $K$-algebra when
 equipped with either of the above topologies and the
operations defined by the equations (\ref{op1}--\ref{op3}). 
\end{prop}

\proof{By the proof of Proposition~\ref{ArtAlgP}, it
suffices to observe that addition and multiplication are
continuous. That is, the sum and product of two Cauchy sequences
should themselves be Cauchy sequences.

Suppose that for $1\leq m\leq j+k$, $(c_{nm}x^{a_{nm}})_{n\geq 0}$
is a Cauchy sequence indexed by $n$.  Let $f_{n}(x)=\sum
_{m=1}^{j}c_{nm}x^{a_{nm}}$, and
$g_{n}(x)=\sum_{m=j+1}^{j+k}c_{nm}x^{a_{nm}}$.  These are
arbitrary Cauchy sequences. Their sum $\sum
_{m=1}^{k+j}c_{nm}x^{a_{nm}}$ is also Cauchy. Their product is a
finite sum of sequences $c_{n,m}c_{n,m'}x^{a_{nm}+a_{nm'}}$. Now,
if either $a_{nm}$ or $a_{nm'}$ decreases (resp.\ increases)
without bound while the other either does the same or is Cauchy,
then their sum decreases without bound. Hence, the sequence is
Cauchy.  Whereas, if both $a_{nm}$ and $a_{nm'}$ are Cauchy, then
so is their sum. Moreover, $c_{nm}$ and $c_{nm'}$ would then both
be Cauchy, so their product is. Hence, the sequence in question is
Cauchy.}

Note that $K[x;R_{1}][y;R_{2}] =
K[y;R_{2}][x;R_{1}]$ as an algebra but {\em not}
as a topological algebra under either topology. 
The expression $K[x,y;R_{1},R_{2}]$ denotes $K[y;R_{2}][x;R_{1}]$ (and
$K[x,y;R]$ denotes $K[y;R][x;R]$).
Expressions like $K[x_{1},\ldots ,x_{n};R_{1},\ldots ,R_{n}]$ (and
$K[x_{1},\ldots ,x_{n};R]$)
are defined in a similar manner. Finally, we define
$$ K[x_{1},x_{2},\ldots;R_{1},R_{2},\ldots  ] = \bigcup_{n\geq 0}K[x_{1},\ldots
,x_{n};R_{1},\ldots ,R_{n}] $$ 
and
$$ K[x_{1},x_{2},\ldots ;R] = \bigcup_{n\geq 0}K[x_{1},.., x_{n};R]. $$

Now, we can give an alternate definition of the Noetherian (resp. Artinian)
algebra. 

\begin{thm}  In the Noetherian (resp.\ Artinian) topology, the
topological completion of the finite algebra $K[x;R]$ is the
Noetherian (resp.\  Artinian) algebra $\hat{K}[x]_{\hat{R}}$
(resp.\ $\hat{K}[x]^{\hat{R}}$) where $\hat{K}$ and $\hat{R}$
are the  
completions of $K$ and $R$ respectively.$\Box $ \end{thm}

\begin{cor} Suppose that $K$ and $R$ are complete. Then in the Noetherian
(resp. Artinian) topology the
completion of ${K[x;R]}$ is $K(x)_{R}$
(resp.\ $K(x)^{R}$).$\Box $  
\end{cor}

\begin{cor}
In the Noetherian topology:\begin{enumerate}
\item The completion of
$K[x_{1},\ldots ,x_{n};R_{1},\ldots ,R_{n}]$ is
$\hat{K}(x_{1;\hat{R}_{1}},\ldots ,x_{n;\hat{R}_{n}})$,
and
\item    The completion of
$K[x_{1},x_{2},\ldots ;R_{1},R_{2},\ldots ]$ is
$\hat{K}(x_{1;\hat{R}_{1}},x_{2;\hat{R}_{2}},\ldots )$
\end{enumerate}
where $\hat{K}$ and $\hat{R_{i}}$
are the completions of $K$ and $R_{i}$ respectively.$\Box $
\end{cor}

\begin{ex}
The completion of $\left( {\bf Q[i]}\right)
[x;{\bf Q}] $ in the Noetherian (resp. Artinian) topology is $ {\bf
C}(x)_{{\bf R}}$ (resp.\ ${\bf C}(x)^{{\bf R}}$) where
$i$ is the imaginary unit.
\end{ex}

Note that $K(x^{R})$ can only be defined in this manner if
$K$ and $R$ are complete. Since Proposition~\ref{Inverses}
supposes that $R$ is Archimedean, and the real numbers are the only complete
Archimedean field \cite{POAS},  we may maintain the
assumption that $R$ is the field of real numbers. 

\ssect{Complex Numbers Raised to a Real Power}

 We need to define
the exponentiation operation for real exponents. Note that the
expression $a^{t}$ is already well defined when $a$ is a
positive real and 
$t$ is an arbitrary real.  It is given by $a^{t}=e^{t\log a}.$

However, there is no unique way to do this when $a$ is an
arbitrary complex number. For example, $\sqrt{-4}=\pm 2i$. In fact,
given any complex number $z$ and real number $t$ we can define a
whole family of values for ``$z^{t}$'' indexed by integer $n$.

\begin{defn}\label{exponentiation} \bld{Exponentiation of a
Complex Number}\index{Exponentiation} Given a
nonzero complex number $z$. Recall that $z$ can be uniquely
written in the form $ae^{i\theta }$ where $i$ is the imaginary
unit $\sqrt{-1}$, $a$ is a positive real number, and $0\leq \theta
<2\pi $. $a$ is called the {\em modulus}\index{Modulus} of $z$ denoted $|z|$
and $\theta $ is called the {\em argument}\index{Argument} of $z$ denoted
$\arg{z}$. The modulus and argument of 0 are both defined to be 0.

Define the $n\th$ value of $z$ to the $t$ power to be $$ z^{t;n}=
a^{t}e^{it(\theta +2n\pi )}.$$

For $z=0$, we define $z^{t;n}=0$ whenever $t$ is positive, and
leave it undefined when $t$ is nonpositive.  \end{defn}

Notice that for all $z,$ $a,$ $n,$ and $m,$ the exponentials $z^{a;n}$ and
$z^{a;m}$ differ by at most a factor of equal to a root of unity. Thus, we can
say that ``$z^{a}$'' is well defined {\em up to multiplication by
factor of modulus $1$.}

\begin{thm}\bld{Characterization of Exponentiation}\index{Exponentiation}
\label{char} $f$ is a 
topological group homomorphism from the group of real numbers
under addition to the group of nonzero complex number under
multiplication if and only if for some nonzero complex number $z$
and some integer $n$, $f(t)=z^{t;n}.$ \end{thm}

\proof{{\bf ($f(t)=z^{t;n}\rightarrow $ Homomorphism)}
Let $s$ and $t$ be real numbers. Set $a=|z|$, and $\theta = \arg
(z)$.  \begin{eqnarray*} f(s)f(t)&=&z^{t;n}z^{s;n}\\ &=&
a^{t}e^{it(\theta +2n\pi )}a^{s}e^{is(\theta +2n\pi )}\\ &=&
a^{s+t}e^{i(s+t)(\theta +2n\pi )}\\ &=&z^{s+t;n}\\ &=&f(s+t).
\end{eqnarray*}

{\bf ($f(t)=z^{t;n}\rightarrow $ Continuous)} Since $z^{t+\delta
;n}=z^{t;n}z^{\delta  ;n}$, it suffices to show  $$ \lim_{\delta \rightarrow
0}z^{\delta ;n }= 1.$$ However, this is immediately obvious.    

{\bf (Only If)} We first characterize the set of such
{\em homomorphisms} from the rational numbers to the nonzero
complex numbers. Since this set is closely related to an extension of the ring
of integers which is of independent interest we  pause and discuss this
ring briefly. Then we determine which homomorphisms
are continuous. These functions then have a unique
continuation to the real numbers.

Let $f(1)=z=ae^{i\theta}$ with $a\geq 0$ and $0\leq \theta < 2\pi $.

By induction we have, $f(na)=f(a)^{n}$ for all integers $n$
and real numbers $a$. Thus, $f(n/m)^{m}=f(n)=z^{n}$. Hence,
$f(n/m)$ is one of the $m$ distinct $m\th$ roots of $z^{n}$.
Thus, $$f(n/m)=a^{n/m}e^{ni(\theta +2k_{n,m}\pi )/m}=z^{n/m;k}$$
for some $k_{nm}$. However, in general this $k_{n,m}$ might depend on
$n$ and $m$.

Now,
\begin{eqnarray*}
f(n/m)&=&f(1/m)^{n}\\
&=&\left(z^{\frac{1}{m};k_{1,m}} \right)^{n}\\
&=&z^{\frac{n}{m};k_{1,m}}
\end{eqnarray*}
by the first part of this proof. Hence, we may assume
without loss of generality that $k_{n_{1},m}=k_{n_{2},m}$ for
all integers $n_{1}$, $n_{2},$ and $m$. We therefore
denote $k_{n,m}$ by $k_{m}$.

Next, observe that the only the value of $k_{m}$ modulo $m$ is
relevant, so we assume that $0 \leq k_{m}<m $.

Finally, note that 
\begin{eqnarray*}
z^{1/m;k_{m}}&=& f(1/m)\\
&=& f(1/nm)^{n}\\
&=& z^{n/m; k_{(nm)}},
\end{eqnarray*}
so we must have $k_{(nm)}\equiv k_{m}$ modulo $n$. Let ${\cal
Z}$ be the set of such sequences.
$$ {\cal Z}=\left\{(k_{n})_{n>0}: 0\leq k_{n}<n\mbox{,
and } k_{nm}\equiv k_{m}\mbox{ modulo $n$}\right\} $$

\digress{We now {\em digress} to consider the algebraic properties of
${\cal Z}$ before continuing on to complete the proof.

${\cal Z}$ is actually a ring since it is the inverse
limit (along the lattice of integers ordered by divisibility)
of the rings ${\bf Z}_{n}$ with respect to the projections:
$$ \begin{array}{rlll}
\phi _{nm,n}: &{\bf Z}_{nm}&\rightarrow & {\bf Z}_{n}\\
&[k]_{nm}&\mapsto &[k]_{n}
\end{array} $$
where $[k]_{n}=k+n{\bf Z}$ is the equivalence class of $k$
modulo $n$.

Note that ${\bf Z}$ can be embedded in in ${\cal Z}$ by
representing $k\in {\bf Z}$ by the sequence $(k_{n})_{n\in
{\bf P}}$ where $k_{n}\equiv k$ modulo $n$.

However, ${\bf Z}$ and ${\cal Z}$ are not identified by this
embedding. For example, if we let $k_{n} \equiv \sum
_{j=1}^{n}j! $ modulo $n$ then $(k_{n})_{n\geq 1}\in {\cal Z}$ but
$(k_{n})_{n\geq 1}$ does not correspond to the sequence of remainders of any
integer.

${\cal Z}$ can be described as a set of
``pseudointegers''\index{Pseudointegers} which 
are characterized by their remainders when divided by
integers. The only requirement on these remainders is that
they be pairwise consistent.}

To achieve any further results, we now must assume that $f$
is continuous. Thus,
\begin{eqnarray*}
1 &=& f(0)\\
&=& \lim_{m\rightarrow +\infty }f(1/m)\\
&=& \lim_{m\rightarrow +\infty }z^{1/m;k_{m}}\\
&=& \left(\lim_{m\rightarrow +\infty }a^{1/m} \right)
\exp \left( \lim_{m\rightarrow +\infty }i\theta /m \right) 
\lim_{m\rightarrow +\infty } \exp \left(2k_{m}i\pi /m \right) \\
&=& \lim_{m\rightarrow +\infty }\exp \left(2k_{m}i\pi /m \right).
\end{eqnarray*}
Hence, the only limit points of the sequence $k_{m}/m$ are 0
and 1. Without loss of generality, 
$$\lim_{m\rightarrow +\infty }{k_{m}/m}=0.$$
Otherwise, eliminate terms from the sequence or multiply by
the image of $-1$ in ${\cal Z}$.

Thus, for all positive integers $j$, $k_{m}/m$ is eventually
less than $1/j$. In particular, $k_{m}$ is eventually equal
to $k_{jm}$. In fact, $k_{m}$ is eventually constant. Thus,
the sequence $k_{m}$ is the representation of some positive
integer $k$, and hence $f(t)=z^{t;k}.$}

\begin{prop}
\label{multrule}
Let $z_{1}$ and $z_{2}$ be complex numbers. Let $t$ be a
real number and let $n$ and $m$ be integers. Then
$$ (z_{1}z_{2})^{t;k}=z_{1}^{t;n}z_{2}^{t;m} $$
where 
$$ k=\left\{ \begin{array}{ll}
n+m&\mbox{if $\arg z_{1}+\arg z_{2}<2\pi $, and}\\
\\
n+m+1&\mbox{if $\arg z_{1}+\arg z_{2}\geq 2\pi $.}
\end{array} \right.$$
\end{prop}

\proo{Let $z_{1}=ae^{\theta i}$ and $z_{2}=be^{\phi
i}$ with $0\leq a,b<2\pi $ and $a,b\geq 0$. Then
\begin{eqnarray*}
z_{1}^{t;i}&=& a^{t}e^{ti(\theta +2n\pi )}\\
z_{2}^{t;j}&=& b^{t}e^{ti(\phi +2m\pi )}\\
(z_{1}z_{2})^{t;k}&=& \left\{ \begin{array}{ll}
(ab)^{t}e^{ti(\theta +\phi +2k\pi )}&\mbox{if $\theta +\phi
<2\pi $, and}\\
\\
(ab)^{t}e^{ti(\theta +\phi +2(k-1)\pi )}&\mbox{if $\theta +\phi
\geq 2\pi .\Box $}
\end{array} \right.
\end{eqnarray*}}

\begin{cor}
The map $z\mapsto z^{t;n}$ is  continuous at all points $z$
which are not nonnegative real numbers.
\end{cor}

\proo{For $z_{0}\not\in {\bf R}^{+}$, $\arg z_{0}\neq
0$, so for $h\in {\bf C}$ near 1, $\arg (z_{0}h)<2\pi $. 
The result now follows from Proposition~\ref{multrule} since
\begin{eqnarray*}
\lim_{z\rightarrow z_{0}}z^{t;n}&=& \lim_{h\rightarrow 1}(z_{0}h)^{t;n}\\
&=& \lim_{h\rightarrow 1}z_{0}^{t;n}h^{t;0}\\
&=&z_{0}^{t;n}.\Box 
\end{eqnarray*}}

\begin{prop}
\label{hyper}
Let $z$ be a complex number. Let  $s$ and $t$ be real
numbers and let $n$ be an integer. Then
$$ z^{st;n}=\left(z^{s;n} \right)^{t;n+k} $$
where $k$ is the greatest integer less than or equal to
${s\arg z}/{2\pi} $.
\end{prop}

\proof{Let $z=ae^{i\theta }$ where $0\leq \theta <2\pi
$, and $a\geq 0$. Then 
\begin{eqnarray*}
z^{st;n}&=&a^{st}e^{sti(\theta +2n\pi )}\\
\left(z^{s;t} \right)^{t;j} &=& a^{st}e^{ti(\phi +2j\pi )}
\end{eqnarray*}
where $0\leq \phi <2\pi $ and $\phi \equiv s(\theta
+2n\pi)$ mod $2\pi $.}

\ssect{Artinian Series to a Real Power}
\label{artrealpow}

In this section, we define the exponentiation of any
nonzero Artinian or Noetherian series whose coefficients
possess a well defined exponentiation. 

As we have seen, these could be series whose coefficients are chosen from
some subset of the complex numbers via
Definition~\ref{exponentiation}. Alternatively, the
coefficients might themselves be Artinian or Noetherian
series of some sort whose exponentiation is defined below. 

First, we must define a few preliminaries:
A {\em multiset}\index{Multiset} $M$ (on a
set $S$) is merely a function (from the set $S$) to the
nonnegative integers. The multiset is {\em
finite}\index{Finite Multiset} if the sum of
its values is finite. This sum is denoted $|M|$.
We use the notation $\sum _{a\in M}f(a)$ to denote
$\sum_{a\in S}M(a)f(a)$. In other words, this is a sum with
multiplicities. Similarly, we calculate products with multiplicities.
$ \prod_{a\in M}f(a)=\prod _{a\in S}f(a)^{M(a)}.$

Next, define the {\em multinomial coefficient}\index{Multinomial Coefficient}
${x \choose M}$ by 
$$ {x\choose M }=\frac{x(x-1)\cdots (x+1-|M|)}{\prod _{j\in
S}M(j)!}.$$

Finally, define the {\em argument}\index{Argument} of an Artinian
or Noetherian series  to be the argument of its leading
coefficient.  

\begin{defn} \bld{Exponentiation of a Series}\index{Exponentiation}
Let $g(x)$ be a
nonzero Artinian (resp.\ Noetherian) series  
whose non-zero coefficients have well defined exponentials.
For example, 
suppose that 
$g(x)$ is a series with complex
coefficients, or that
$g(x)$ is a series whose coefficients are
series with complex coefficients, or that
$g(x)$ is a series whose coefficients are
series whose coefficients are series with complex coefficients and so on.

Let $d=\deg(g(x))$, and  $c_{a}=[x^{a}]g(x)$. Then for
all integers $n$ and real numbers $t$ define the $g(x)$ to the power $t$
indexed by $n$ by the sum
\begin{equation}\label{ftoa}
g(x)^{t;n}=c_{d}^{t;n} x^{dt} \sum _{M} {t\choose M}
\left(\prod _{a\in M} \frac{c_{a+d}}{c_{d}}x^{a} \right)
\end{equation}
over all multisets $M$ over the set ${\bf
R}-\{0 \}$ of nonzero real numbers.

For the exponentiation of $g(x)=0$, we follow the same
convention as in Definition~\ref{exponentiation}.
\end{defn}

\begin{prop}
\label{welldef-exp}
Let $g(x)$ be a nonzero Artinian (resp.\ Noetherian) series.
Let $t$ be a real number, and $n$ an integer. Then
$g(x)^{t;n}$ is a well defined Artinian (resp.\ Noetherian)
series.
\end{prop}

\proof{It suffices to show
that $[x^{b}](g(x)^{t;n})$ is well defined for all real
numbers $b$, and that there are finitely many $a\leq b$
(resp.\ $a\geq b$) such that $[x^{a}](g(x)^{t;n})$
nonzero.

Now, there are finitely many $a\leq b-dt$ (resp.\ $a\geq
b-dt$) such that $c_{a}\neq 0$. Denote them $a_{1},\ldots
,a_{k}$. 

We claim that the summation yields only finitely many terms
of degree less (resp.\ greater) than or equal to $b$. Thus,
 $[x^{b}](g(x)^{t;n})$ is well defined, and there are
only finitely many $a\leq b$ (resp.\ $a\geq b$) such that
 $[x^{a}](g(x)^{t;n})$ is nonzero.

The claim is true since for the only summands which
contribute terms of interest are those indexed by a multiset
$M$ such that\begin{enumerate}
\item When $a\in \{a_{1},\ldots ,a_{k} \} , $ we have $M(a)=0$.
\item Conversely, $M(a_{j})\leq \frac{b-dt}{a_{j}}.$
\end{enumerate}
Obviously, there are only finitely many such multisets.}

Let $t$ be a real number, and let $n$ and $m$ be integers.
Then $g(x)^{t;n}$ and $g(x)^{t;m}$ differ only by a factor
which is a root of unity. Thus ``$g(x)^{t}$'' is well
defined up to multiplication by complex numbers of modulus
one.

In particular, $g(x)^{-1;n}=g(x)^{-1;m}$ for all $m$ and
$n$. It will be shown (Proposition~\ref{Charact2}) that the
exponentiation $g(x)^{-1;n}$ gives an explicit formulation
of the reciprocal of $g(x)$ which had been shown to exist
via Proposition~\ref{Inverses}.

\begin{prop}
\label{Duality2}
The isomorphism $\iota $ defined in
Proposition~\ref{Duality} preserves exponentiation. That is, informally,
$$ g^{t;n}(1/x)=g(1/x)^{t;n} $$
for all Noetherian or Artinian series $g(x)$, and all
integers $n$ and real numbers $t.\Box $
\end{prop}

\begin{thm}\bld{Characterization of Exponentiation}\index{Exponentiation}
\label{Charact2} 
 $\sigma $ is a topological group homomorphism from the group
of real numbers under addition to the group of nonzero
Artinian (resp.\ Noetherian) series under 
multiplication if and only if for some nonzero series
$g(x)$ and some integer $n$,  $\sigma (t)=g(x)^{t;n}.$
\end{thm}

\proof{{\bf (If)} The function
$c_{d}^{t;n}$ is continuous by Theorem~\ref{char}, and 
polynomials are continuous. So for this direction of need only show
$$ f(x)^{a+b;n}=\left(f(x)^{a;n} \right)\left(f(x)^{b;n} \right). $$
This is true for $f(x)$ a monomial by Theorem~\ref{char}, so it 
suffices to show for $f(x)=1+\sum_{a>0}c_{a}x^{a}$.
\begin{eqnarray*}
f(x)^{a;0}f(x)^{b;0}&=&
\left(\sum _{M} {a\choose M }\prod_{a\in M}c_{a}x^{a}\right)
\left(\sum _{N} {b\choose P }\prod_{b\in M}c_{b}x^{b}\right)\\
&=& \sum _{P}\left(\sum _{M+N=P}{a\choose M }{b\choose M}
\right)\prod _{a\in P}c_{a}x^{a},
\end{eqnarray*}
and this is equal to 
$$ \sum _{P}{a+b\choose P}\sum _{a\in P}c_{a}x^{a} $$ 
 by the analytic version of the multivariate
Van der Monde convolution.

{\bf (Only If)} Up to a constant, every series has exactly
one $n\th$ root.  Since there are $n$ choices of constant by
Theorem~\ref{char},
there is a choice of $k$ such that
$\sigma (a)=\sigma (1)^{a;k}$ for all rational numbers $a$.
By continuity, this holds for all real numbers $a.$}

\begin{prop}
\label{Mult2}
Let $f(x)$ and $g(x)$ be Artinian (resp.\ Noetherian) series.
Then for
all integers $j$ and $k$ and real numbers $t$,
$$ (f(x)g(x))^{t;n}=\left( f(x)^{t;k} \right)
\left(g(x)^{t;j} \right) $$
where
$$ n=\left\{ \begin{array}{ll}
j+k&\mbox{if $\arg f(x)+\arg g(x)<2\pi $, and}\\
\\
j+k+1&\mbox{if $\arg f(x)+\arg g(x)\geq 2\pi $.}
\end{array} \right. $$
\end{prop}

\proof{By Proposition~\ref{multrule}, this holds for
constants, so we may assume without loss of generality that
$f(x)$ and $g(x)$ have leading term 1. Let
$c_{a}=[x^{a}]f(x)$, $d_{b}=[x^{b}]g(x)$, and
$e_{b}=\sum _{a_{1}+a_{2}=b}c_{a_{1}}d_{a_{2}},$ so
that $f(x)g(x)=1+\sum _{b>0}e_{b}x^{b}$ (resp.
$f(x)g(x)=1+\sum _{b<0}e_{b}x^{b}$).

Now, \begin{eqnarray*}
\sum _{P}{t\choose P }\sum _{b\in P}e_{b}x^{b}&=& 
\sum _{P}{t\choose P }\sum _{b\in P}\sum _{a_{1}+a_{2}=b}c_{a_{1}}d_{a_{2}}x^{b}\\
&=&\sum _{P}{t\choose P }\sum_{M+N=P} 
\left(\prod_{a}{P(a)\choose M(a)} \right) 
\left(\prod _{a_{1}\in M} c_{a_{1}}x^{a_{1}}\right)
\left(\prod _{a_{2}\in N} c_{a_{2}}x^{a_{2}}\right)\\
&=&\sum_{M,N}
\frac{t(t-1)\cdots(t+1-|M|-|N|)}{\prod_{a}M(a)!N(a)!}
\prod _{a_{1}\in M}c_{a_{1}}x^{a_{1}} \prod _{a_{2}\in N}d_{a_{2}x^{a_{2}}}\\
&=& f(x)^{t;i}g(x)^{t;j}
\end{eqnarray*}
since $t(t-1)\cdots (t+1-m)$\index{Lower Factorial} is a sequence of
polynomials of binomial type.\cite{RR}}

The following Lemma demonstrates that \eref{ftoa} is much more general
than previously indicated.

\begin{lem}\label{cornish}
Let $f(x)$ be an Artinian (resp.\ Noetherian) series of
degree $d$. Suppose $f(x)$ is given by the following 
convergent expansion
$$f(x)=\sum_{j\geq 0}c_{j}x^{a_{j}}$$
where $(a_{j})_{j\geq 0}$ is a sequence of not necessarily
distinct  real numbers such that $a_{0}=d$ and for $j$
positive, $a_{j}>d$ (resp.\ $a_{j}<d$). Then $f(x)^{t;n}$ is given by the sum
$$ f(x)^{t;n}=(c_{d})^{t;n} x^{a_{0}t} \sum _{M} {t\choose
M} \prod _{j\in M} \frac{c_{j}}{c_{0}}x^{a_{j}-a_{0}}$$
 over all finite multisets of positive integers. 
\end{lem}

\proof{Suppose that $a_{i},a_{j},a_{k},\ldots $ are
all equal to $a$. Then any selection of $i,j,k,\ldots $ for
the multiset $M$ contributes a factor of $x^{a}$ just as a
selection of $a$ normally would. By the multinomial theorem,
$$ \sum _{|M|=n}{n\choose M }\prod _{l \in M}
\frac{c_{l}}{c_{0}} =
\left(\frac{c_{i}+c_{j}+c_{k}+\cdots }
{c_{0}}\right)^{n} $$ 
where $n$ is a nonnegative integer and the sum ranges over
all $n$-element multiset over $\{i,j,k,\ldots  \}$.
Now, multiply both sides by $t(t-1)..(t+1-n)/n!.$}

\begin{prop}
\label{othercts}
Let $t$ be a real number and $n$ an integer. Then the map
$f(x)\mapsto f(x)^{t;n}$ is continuous at all
nonzero Artinian (resp.\ Noetherian) series $f(x)$
such that $\arg f(x)\neq 0.\Box $
\end{prop}

\begin{prop}
\label{iterExp}
Let $f(x)$ and $g(x)$ be an Artinian (resp.\ Noetherian)
series. Then for all integers $n$ and 
real numbers $s$ and $t$
$$ f(x)^{st;n}= \left(f(x)^{s;n} \right)^{t;n+k} $$
where $k$ is the greatest integer less than or equal to
$s\arg f(x)/2\pi$. 
\end{prop}

\proof{This is true for constants by
Proposition~\ref{hyper}, so without loss of generality, the
leading term of $f(x)$ is 1 and $n$ and $k$ are zero. Let
$c_{a}=[x^{a}]f(x)$. By Lemma~\ref{cornish}, 
$$ \left(f(x)^{s;0} \right)^{t;0} = 
\sum_{\cal M}{t\choose {\cal M} }
\prod_{ M\in {\cal M}} {s\choose M} \prod _{a\in M}c_{a}x^{a} $$
where the sum is over finite multisets ${\cal M}$ of finite
nonempty multisets $M$ of nonzero real numbers. We invert
the order of summation:
$$ \left(f(x)^{s;0} \right)^{t;0} = 
\sum _{N}\left(\sum _{{\cal M}}{t\choose {\cal M} } \prod
_{M\in {\cal M}}{s\choose M} \right) \prod _{a\in N}c_{a}x^{a} $$
where the inner sum is over  multisets ${\cal M}$ as above
such that 
$$N=\sum _{M\in {\cal M}}M.$$
It remains now to show that 
$$\sum _{{\cal M}}{t\choose {\cal M} } \prod
_{M\in M}{s\choose M}={st\choose N}.$$
When $s$ and $t$ are nonnegative integers, this holds by
combinatorial reasoning; both sides count the number 
of partitions of an $st$-set whose multiset of block sizes
is $N$. Moreover, since this is an identity of {\em
polynomials} we now know that it holds for all $s$ and $t$.}

\begin{prop}
\label{degree}\label{dt}
Let $f(x)$ be an Artinian (resp.\ Noetherian) series of
degree $d$. Let $t$ be a real number, and $n$ be an integer.
Then $\deg(f(x)^{t;n})=td.\Box $
\end{prop}

\begin{cor}
Let $f(x)$ be an Artinian (resp.\ Noetherian) series whose
degree is nonzero. Let $n$ be an integer. Then the set
$\{f(x)^{t;n}: t\in {\bf R}\}$ is a $K$-pseudobasis for the
Artinian algebra (resp.\ Noetherian algebra).$\Box $
\end{cor}

Actually, much more is true of these pseudobases as we 
see in the next section.

\ssect{Composition of Series}\label{cos}

\begin{defn}\bld{Composition of Series}\label{CompArt}\index{Composition}
Given two Artinian (resp.\ Noetherian)
series $f(x)$ and $g(x)$ such that $g(x)$ has {\em positive} degree, and an
integer $n$. Suppose that
$f(x)=\sum _{a}c_{a}x^{a}$. Then define the $n\th$
{\em composition} of $f(x)$ with $g(x)$ to be 
$$ f(g;n)=\sum _{a}c_{a}g(x)^{a;n}. $$
\end{defn}

Note that $g(x)$ must have positive degree regardless of whether the
series is Artinian or Noetherian. In other words, the isomorphism $\iota $ from
Propositions \ref{Duality} and \ref{Duality2} does not preserve composition.

\begin{thm}
\bld{Characterization of Composition}\label{CharComp}\index{Composition}
Let $g(x)$ be an Artinian (resp. Noetherian) series 
positive degree. Then for all integers $n$,
$f(x)\rightarrow f(g;n)$ is a continuous field automorphism of the
Artinian (resp.\ Noetherian) algebra which fixes all constants.

Conversely, any continuous field automorphism of the
Artinian (resp.\ Noetherian) algebra which fixes all
constants is of the form $f(x)\rightarrow f(g;n)$ for some
Artinian (resp.\ Noetherian) series of positive degree $g(x)$
and some integer $n$.
\end{thm}

\proof{We must show the following for both Artinian
and Noetherian series:
\begin{itemize}
\item $f(g;n)$ is a well defined Artinian (resp.\ Noetherian) series. 
\item Composition is continuous. 
\item Composition preserves addition. 
\item Composition preserves constants. 
\item Composition preserves multiplication. 
\item The Converse. 
\end{itemize}

{\bf (Well Defined)} Suppose $\deg(g(x))=d>0$, and
$f(x)=\sum _{a}c_{a}x^{a}$. Then 
$$[x^{a}](f(g;n))=\sum _{b\leq a/d} c_{b}[x^{a}](g(x)^{b;n})$$
(resp.\ $\sum _{b\geq a/d}c_{b}[x^{a}](g(x)^{b;n})$)
 is a finite sum since $c_{b}$ is an Artinian (resp.
Noetherian) sequence.

The sequence $[x^{a}](f(g;n))$ is itself Artinian (resp.
Noetherian) since there are only finitely many $b\leq a/d$
(resp.\ $b\geq a/d$) such that $c_{b}\neq 0$ and each $g(x)^{b;n}$ contributes
only finitely many terms of order less (resp.\ greater) than $x^{a}$.

{\bf (Continuous)} It is a convergent sum of continuous
functions by Theorem~\ref{Charact2}.

{\bf (Addition)} Obvious.

{\bf (Constants)} If $f(x)=c$, then $f(g;n)=cg(x)^{0;n}=c$.

{\bf (Multiplication)} Immediately follows from
Theorem~\ref{Charact2}.

{\bf (Converse)} Such an automorphism $\sigma $ could be composed with
the group automorphism $\theta :a\mapsto x^{a}$. The resulting map
$\theta \circ \sigma $ satisfies the conditions of Theorem~\ref{Charact2}, so
we must have $\sigma x^{a} = g(x)^{a;n}$ for some series $g(x)$
and integer $n$. By continuity, we know  $g(x)^{a;n}$ tends
tends to zero as $a$ tends to $+\infty $ (resp.\ $-\infty $).
Thus, $\deg(g(x))\geq 0$ and $\sigma f(x)=f(g;n)$.}

\begin{prop}
\label{deg-comp}
For all Artinian (resp.\ Noetherian) series $f(x)$ and
$g(x)$ and integers $n$,
$$\deg(f(g;n))=\deg(f(x))\deg(g(x)).$$
\end{prop}

\proof{Proposition~\ref{dt}.}

As opposed to the composition of functions or of the usual sort of formal power
series, this composition is not necessarily associative.

\begin{ex}
Let $b$, and $c$ be small positive
real numbers. Set $f(x)=x$, $g(x)=x^{b}$, and
$h(x)=x^{c}$. Then 
$$ (f(g;1))(h;1)=x^{bc}e^{2\pi i(1+c(b+1))} $$
and
$$ f(g(h;1);1)=x^{bc}e^{2\pi ic(1+2b)} .$$
Their quotient is $e^{2\pi i(1-b)}$ which is not equal to one.
\end{ex}

However, all is not lost; indeed, composition by is associative
in the case we are most interested in: delta series; for in
\cite{ch4}, the relevant computations involve
compositions with  delta series.

\begin{thm}
Let $m$ be an integer, and let $f(x)$, $g(x)$, and $h(x)$ be
Artinian (resp.\ Noetherian) series such that $g(x)$ and
$h(x)$ are delta series and $\arg g(x)+\arg h(x)<2\pi $. Then
$$ (f(g;0))(h;m)=f(g(h;m);0). $$
\end{thm}

\proof{By linearity and continuity, it suffices to
consider $f(x)=x^{a}$. Suppose $c_{b}=[x^{b}]f(x)$ and
$d_{b}=[x^{b}]g(x)$. Then
\begin{eqnarray*}
(f(g;0))(h;m)&=& (g(x)^{a;0})(h;m)\\
&=& \left(c_{1}^{a;0}x^{a} \sum_{M}{t\choose M}\prod _{b\in
M}\frac{c_{b+1}x^{b}}{c_{1}} \right) (h;m)\\
&=& c_{1}^{a;0}h(x)^{a;m} \sum_{M}{t\choose M}\prod _{b\in
M}\frac{c_{b+1}h(x)^{b;m}}{c_{1}}.
\end{eqnarray*}
Whereas,
\begin{eqnarray*}
f(g(h;m);0)&=& f\left(c_{1}h(x)+\sum _{b>1}c_{b}h(x)^{b;m};0 \right)\\
&=& (c_{1}h(x))^{a ;0} \sum_{M}{t\choose M}\prod _{b\in
M}\frac{c_{b+1}h(x)^{b;m}}{c_{1}}
\end{eqnarray*}
by Lemma~\ref{cornish}. Finally,
$c_{1}^{a;0}h(x)^{a;m}=(c_{1}h(x))^{a;m}$ by
Proposition~\ref{Mult2}.}

\begin{por}
Let $m$ and $n$ be integers, and let $f(x)$, $g(x)$, and $h(x)$ be
nonzero Artinian (resp. Noetherian) series with $\deg(g(x)),\deg(h(x))>0$. Then
the series $ (f(g;0))(h;m)$ and  $f(g(h;m);0)$ only differ by a
factor which is a root of unity.$\Box $
\end{por}

\begin{cor}
The set of Artinian (resp.\ Noetherian) delta series is a
monoid with respect to 0-composition.$\Box $
\end{cor}

\begin{prop}
 For all integers $n$ and Artinian (resp.\ Noetherian)
series $f(x)$, the map $g(x)\mapsto f(g;n)$ is continuous at
all  series $g(x)$ of positive degree such that $\arg
g(x)\neq 0$.
\end{prop}

\proof{Proposition~\ref{othercts}.}

\section{Derivative}

\begin{defn} \bold{Derivative} \label{dsec}
We  define the {\em derivative}
(with respect to $x$) to be the continuous linear map on the
Artinian (resp.\ Noetherian) algebra denoted $D$ or $\frac{d}{dx}$
such that $$ Dx^{a}=ax^{a-1}. $$ The derivative of $f(x)$ is
denoted $f'(x).$ \end{defn}

\begin{prop}\bold{Derivation} \label{derivP1}
The derivative is a derivation. That is,
$$ D(f(x)g(x))=f'(x)g(x)+f(x)g'(x).$$ \end{prop}

\proo{Let $c_{a}=[x^{a}]f(x)$ and
$d_{a}=[x^{a}]g(x)$. Then 
\begin{eqnarray*} 
f'(x)g(x)+f(x)g'(x)&=&\sum_{b}\sum
_{a_{1}+a_{2}=b}(a_{1}+a_{2}) c_{a_{1}}d_{a_{2}} x^{b-1}\\ 
&=& D(f(x)g(x)).\Box
 \end{eqnarray*}}

\begin{prop}
\label{der-pow}
Let $f(x)$ be a nonzero Artinian (resp.\ Noetherian) series. Let $t$
be a real number and $n$ be an integer. Then
$$ D(f(x)^{t;n}) = tf'(x)f(x)^{t-1;n} $$
\end{prop}

{\em Proof:} Let $f(x)=\sum_{b}c_{b}x^{b}$. Suppose $f(x)$
has degree $d$. Then
\begin{eqnarray*}
tf'(x)f(x)^{t-1;n}&=& 
t \left( \sum _{a} ac_{a}x^{a-1} \right)
\left(c_{d}^{t-1;n}x^{d(t-1)} \right)
\sum _{M}{t-1\choose M }\prod _{a\in M}\frac{c_{a+d}}{c_{d}}x^{a}\\
&=& c_{d}^{t;n}
\left( dtx^{dt-1} +tx^{dt}\sum _{a\neq 0} \frac{(a+d)c_{a+d}}{c_{d}}x^{a-1}
\right) 
\sum _{M}{t-1\choose M }\prod _{a\in M}\frac{c_{a+d}}{c_{d}}x^{a}\\
&=&  (\D c_{d}^{t;n}x^{dt}) \left(\sum _{M} {t\choose M}
\prod_{a\in M}\frac{c_{a+d}}{c_{d}}x^{a}  \right)
 +c_{d}^{t;n}x^{dt}\D \left(\sum _{M} {t\choose M } \prod_{a\in
M}\frac{c_{a+d}}{c_{d}}x^{a} \right)\\
&=& Df(x)^{a;n}.\Box 
\end{eqnarray*}

\begin{prop}\bold{Chain Rule}\index{Chain Rule}
\label{(Chain Rule)}
For all valid $n$-compositions of series $f(x)$ and $g(x)$,
$$ \D (f(g;n))=f'(g;n)g'(x). $$
\end{prop}

\proof{By linearity and continuity, it suffices to
consider the case $f(x)=x^{a}$. However, this case was
treated by Proposition~\ref{der-pow}.}

\section{LaGrange Inversion}\label{AppLag}
In this section, we show the existence of compositional
inverses to nonzero Artinian (resp.\ Noetherian) series of
positive degree, and we give an explicit formula for
their coefficients.

\begin{prop}
\label{semi}
 Let $f(x)$ be a  Artinian (resp.\ Noetherian) series whose
leading coefficient is a positive real number, and whose
degree is positive. Then $f(x)$ possesses a left (resp.\ right)
$0$-compositional inverse 
denoted $f^{(-1;n)}(x)$ for all integers $n$. 

Thus, the set of such series is a semigroup under $0$-composition.
\end{prop}

\proof{We construct $f^{(-1;n)}(x)$ from its coefficients
$c_{a}=[x^{a}]f^{(-1;n)}(x)$. 

Let $b=\deg(f(x))$. For $a<1/b$, set $c_{a}=0$. Then
set $c_{1/b}=d_{b}^{(1/b);0}$. Next, recursively
define $c_{a}=-z_{a}^{(1/b);0}$ where
$z_{a}=[x^{ab}]f(\sum_{e<a}c_{e}x^{e};n).$

As in the proof of Proposition~\ref{Inverses}, $c_{a}$ is an Artinian series
since ${\bf R}$ is Archimedean. 

Let $g(x)=\sum _{a}c_{a}x^{a}$. We now immediately have $x=f(g;n)$.}

\begin{thm}
\bold{LaGrange Inversion} \label{LaG}
Let $f(x)$ be a delta series with leading coefficient a
positive real number. Then
\begin{equation}\label{lageq}
 [x^{a}] \left(f^{(-1;0)}(x)
\right)^{b;0}=b[x^{a-b}]\left(\frac{x}{f(x)} \right)^{a;0}
\end{equation}
for all real numbers $a$ and $b$.
\end{thm}

Note  that for $b=1$, \eref{lageq} immediately gives the
coefficients of $f^{(-1;0)}(x)$.

\proo{First, observe that $[x^{-1}]\D  g(x)=0$ for
any series $g(x)$. 

Now, let $d_{a}= [x^{a}] \left(f^{(-1;0)}(x)
\right)^{b;0}$. Then
\begin{eqnarray*}
x^{b}&=&\sum_{c\geq b}d_{c}f(x)^{c;0}\\
bx^{b-1}&=&\sum_{c\geq b}cd_{c}f(x)^{c-1;0}f'(x)\\
bx^{b-1}/f(x)^{a;n} &=& \sum_{c\geq b}cd_{c}f(x)^{c-a-1;0}f'(x).
\end{eqnarray*}
Now, for $a\neq c$, $f(x)^{c-a-1;0}f'(x)$ is the derivative
of $\frac{1}{c-a}\D f(x)^{c-a-1;0}$. Hence, by the above observation
\begin{eqnarray*}
[x^{-1}] \frac{kx^{k-1}}{f(x)^{a;0}} &=& [x^{-1}]
ad_{a}\frac{\D f(x)}{f(x)}\\
&=& ad_{a}\\
&=& a[x^{a}]\left(f^{(-1;0)}(x) \right)^{b;0}.\Box 
\end{eqnarray*}}

\section{Symmetric Artinian and Noetherian Series}
\sssect{The Symmetric Algebras}

As a final direct application of theory of Artinian (resp.~Noetherian) series,
we now demonstrate that symmetric Artinian (resp.~Noetherian)
series can be studied in the same way as symmetric functions with integral
exponents. 

A multivariate Artinian (resp.~Noetherian) series is {\em
symmetric}\index{Symmetric Series} if it is invariant under exchanges of
variables. 
Let the $\Lambda (n)$ denote the symmetric part of the multivariate Noetherian
algebra $K(x_{1},x_{2},\ldots ,x_{n})_{{\bf R}^{+}}$, and let $\Lambda ^{*}(n)$
denote the symmetric part of the multivariate Artinian algebra
$K(x_{1},x_{2},\ldots ,x_{n})^{{\bf R}^{+}}$ where ${\bf R}^{+}$ is the set of
nonnegative real numbers. 

Now, let $\Lambda _{a}(n)$ and $\Lambda _{a}^{*}(n)$ consist of those series in
$\Lambda $ and $\Lambda ^{*}$ respectively which are homogeneous of degree $a$.
Thus, $\Lambda (n)$ and $\Lambda ^{*}(n)$ are graded by the modules $\Lambda
_{a}(n)$ and $\Lambda _{a}^{*}(n)$. 

By Porism~\ref{poor}, a symmetric series is invertible if and only if its
constant term is zero.

Also note that given any Artinian (resp.~Noetherian) symmetric series
$p(x_{1},x_{2},\ldots ,x_{n})$ without constant term, any Artinian
(resp.~Noetherian) series $q(y)$, and any integer $k$; the composition
$ q(p(x_{1},x_{2},\ldots ,x_{n});k) $ is well defined.

Let us say that  the ring $K$ is an {\em exponential ring}\index{Exponential
Ring} over the monoid $R$ if it is associated  with a family of ring
isomorphisms indexed by pairs $(a;k)$ where $a\in R$, and $k$ is an integer;
the image of a point $x\in K$ 
under the morphism indexed by $(a;k)$ is $x^{a;k}$. We further require that
$x^{a;k}x^{b;k}= x^{a+b;k}$. If $K$ is endowed with a topology, an algebra
structure over another exponential ring, or a grading, then these maps
must be continuous with respect to that topology, respect the scalar
multiplication, and preserve the grading.

Above we first showed that the complex numbers formed an exponential ring over
the real numbers. Next, we showed that if $K$ was an exponential ring over the
real numbers, then so is $K(x)^{{\bf R}}$ and $K(x)_{{\bf R}}$. This implies
that the corresponding multivariate and even infinite multivariate Noetherian
and Artinian algebras are exponential rings over the real numbers.
Clearly, by the above remarks the algebras of Artinian and Noetherian symmetric
series are graded exponential algebras over the nonnegative real numbers.

Before continuing to consider symmetric series, we pose the following question:
\begin{prob}
Is it possible to work with
Artinian and Noetherian transcendence bases in the same manner as one works
with algebraic transcendence bases?
\end{prob}

Now, the map $\Lambda(n) \rightarrow \Lambda(n-1)$ (resp.~$\Lambda(n)
\rightarrow \Lambda(n-1)$) formed by setting the $n\th$ variable equal to zero
is obviously an exponential ring homomorphism. (That is, it is an algebra
homomorphism and it respects the exponential maps.) Thus, we can proceed to
invoke inverse limits along the category of graded exponential algebras and
thus define the Artinian and Noetherian symmetric algebras in infinitely many
variables. 
\begin{eqnarray*}
\Lambda &=& \lim_{{\displaystyle \leftarrow \atop n}} \Lambda (n)\\*
\Lambda^{*} &=& \lim_{{\displaystyle \leftarrow \atop n}} \Lambda^{*} (n).
\end{eqnarray*}

Note that this inverse limit must be taken as a {\em graded} exponential
algebra. Otherwise, one would unintentionally include 
 series of arbitrarily large degree such as the infinite product
$\prod_{i=0}^{\infty } (1+x_{i})$.

We may not define a symmetric algebra which is an exponential ring over all of
${\bf Z}$---at least not in an infinite number of variables. Consider the
product of the symmetric series $\sum_{i} x_{i}$ and $\sum _{i} x_{i}^{-1}$.
Their product in a finite number of variables $n$ is $n+2\sum_{1\leq i<j\leq
n}x_{i}x_{j} $. Thus, the product is not well behaved under changes in the
number of variables. In particular, the product is not well defined when the
set of variables is infinite.

\sssect{The Monomial Symmetric Series}

Normally, the monomial symmetric function is indexed by an integer partition
$\lambda $; that is, a nonincreasing vector with finite support. (The nonzero
entries of $\lambda $ are called its parts, and $\ell(\lambda )$ is the number
of parts of $\lambda $).  Here this is not appropriate, so we define a {\em
real partition}\index{Partition, Real} to be  a  vector $\beta $ of nonnegative
real numbers such that $\beta $ is nonincreasing and has finite support. 

\begin{defn}\label{monomial}\bold{Monomial Symmetric Series}
For each real partition $\beta$, define the monomial symmetric function by the
sum 
$$ m_{\beta}({\bf x})= \sum_{\alpha }{\bf x}^{\alpha } $$
over all {\em distinct} permutations $\alpha $ of the partition $\beta $ where
${\bf x}=\{x_{1},x_{2},\cdots  \}$ and ${\bf x}^{\alpha }=x_{1}^{\alpha _{1}}
x_{2}^{\alpha _{2}}\cdots $.
\end{defn}

Clearly, the monomial symmetric series are well defined, they are indeed
symmetric, and they generalize the classical monomial symmetric functions.
Moreover, the $m_{\beta}({\bf x})$ where the $\beta_{i}$ sum to $a$
form a basis for the $\Lambda _{a}$ (resp.~$\Lambda^{*}_{a}$).
Thus, the collection of all $m_{\beta}({\bf x})$ form a basis for $\Lambda$
(resp.~$\Lambda^{*}$). 

\sssect{Transcendence Bases}

Recall the definition of the {\em elementary symmetric
function} 
$e_{n}({\bf x}) $ and {\em complete symmetric function} $h_{n}({\bf
x}) $. They are defined explicitly by
the sums
\begin{eqnarray*}
e_{n}({\bf x}) &=& \sum _{\scriptstyle \mu \in
{\cal P}^{*} \atop \scriptstyle \ell (\mu )=n} m_{\mu}({\bf y})\\
h_{n}({\bf x}) &=& \sum _{\scriptstyle \lambda \in
{\cal P} \atop \scriptstyle \ell (\lambda )=n} m_{\lambda}({\bf y})
\end{eqnarray*}
 over integer partitions $\mu$ with
distinct parts and over all integer partitions $\lambda $,
and they are defined implicitly by the generating functions
\begin{eqnarray*}
\prod _{n\geq 1}(1+x_{n}y)&=&\sum _{n\geq
0}e_{n}({\bf x})y^{n}\\
\prod _{n\geq 1}(1-x_{n}y)^{-1}=\sum _{n\geq
0}h_{n}({\bf x})y^{n}.
\end{eqnarray*}
Now, generalize to Artinian and Noetherian series.
\begin{defn}\bold{Elementary and Complete Symmetric Series} 
Let $\beta$ be a real partition. Define the {\em elementary} or {\em complete}
symmetric series by the product
\begin{eqnarray*}
e_{\beta}({\bf x}) &=& e_{1}({\bf x})^{\beta_{2}-\beta_{1}} e_{2}({\bf
x})^{\beta_{3}-\beta_{2}}\cdots \\*
h_{\beta}({\bf x}) &=& h_{1}({\bf x})^{\beta_{2}-\beta_{1}} h_{2}({\bf
x})^{\beta_{3}-\beta_{2}}\cdots .
\end{eqnarray*}
\end{defn}

These series are well defined and symmetric, and is a generalization of the
definition of $e_{\lambda ^{*}}({\bf x})$ and $h_{\lambda ^{*}}({\bf x})$.
Each series can be expressed as a sum 
$$m_{\beta}({\bf x}) +\sum_{\gamma >\beta}a_{\gamma \beta}m_{\gamma}({\bf x})$$
over real partitions $\gamma $ occurring after $\beta$ in reverse
lexicographical order. Thus, the $e_{\beta}({\bf x})$ and the $h_{\beta }({\bf
x})$ are bases for $\Lambda $ and $\Lambda ^{*}$.

In general, we say that a subset $S$ of an exponential ring $K$ over the
monoid $R$ is an {\em Artinian} (resp.~{\em Noetherian}) {transcendence
basis}\index{Transcendence Basis} if for a given $n$, every element $x\in R$ is
represented by a unique Artinian (resp.~Noetherian) series over the elements of
$S$. 

By the above reasoning, the elementary and complete symmetric functions
$e_{n}({\bf x}) $ and $h_{n}({\bf x})$ each form an Artinian transcendence
basis for $\Lambda ^{*}$, and a Noetherian transcendence basis for $\Lambda $.
Thus, the classical involution of the algebra of symmetric functions $\omega :
h_{n}({\bf x})\leftrightarrow e_{n}({\bf x})$ extends to a well defined
involution of the exponential algebras of symmetric series $\Lambda $ or
$\Lambda ^{*}$.

Thus, we can define the {\em forgotten symmetric series}\index{Forgotten
Symmetric Series} via the identity $f_{\beta }({\bf x})=\omega m_{\beta }({\bf
x})$.

Define the {\em power sum symmetric function} $p_{n}({\bf x})=\sum _{i}
x_{i}^{n}$. As in \cite{Mac}, we can uniquely express the $p_{n}({\bf x})$ in
terms of the  
$e_{n}({\bf x})$ and visa-versa. Thus, the power sum symmetric functions
constitute yet another Noetherian transcendence basis for $\Lambda $ and
Artinian transcendence basis for $\Lambda ^{*}$.

Thus, the power sum symmetric series which are products of the power sum
symmetric functions
$$ p_{\beta }({\bf x})= p_{1}({\bf x})^{\beta _{1}- \beta _{2}} p_{1}({\bf
x})^{\beta _{1}- \beta _{2}}  \ldots $$
form a basis for $\Lambda $ and $\Lambda ^{*}$. 
To compute the action of $\omega $ with respect to this basis, note that
$\omega p_{n}({\bf x})=(-1)^{n-1}p_{n}({\bf x})$. Thus,
$$ \omega p_{\beta }({\bf x})= (-1)^{\beta _{2}+\beta _{3}+\cdots }p_{\beta
}({\bf x}).$$

\end{document}